\documentclass[12pt,a4paper, reqno]{amsart}
\usepackage{amsmath,amsfonts,amssymb,mathrsfs}
\usepackage{mathtools}
\usepackage{graphicx}
\usepackage{comment}
\usepackage{xfrac}
\usepackage{subfigure}
\usepackage{bbm}
\usepackage{epstopdf}
\usepackage{pstricks}
\usepackage{pst-node}
\usepackage[linktocpage=true,colorlinks=true,linkcolor=blue,citecolor=red,urlcolor=magenta,pdfborder={0 0 0}]{hyperref}
\usepackage{fancyhdr}
\usepackage{mathscinet}
\usepackage{dutchcal}
\usepackage{xcolor}
\usepackage[english]{babel}
\usepackage{mathscinet}

\usepackage{chngcntr}

\usepackage{appendix}

\definecolor{forestgreen}{rgb}{0.13, 0.55, 0.13}
\definecolor{anna}{rgb}{0.01, 0.28, 1.0}

\newtheorem{theorem}{\bf Theorem}[section]

\newcommand{\RR}{\mathbb{R}}

\newcommand{\OO}{\mathbb{O}}
\newcommand{\I}{\mathbb{I}}

\def \L {\mathscr{L}}
\def \K {\mathscr{K}}

\def \o {{\omega}}

\def \b {{\beta}}

\def \g {{\gamma}}

\def \epsilon {{\varepsilon}}

\def \k {{\kappa}}

\def \s {{\sigma}}

\def \phi {{\varphi}}

\def\p{\partial}

\usepackage{mathtools}
\DeclarePairedDelimiter{\abs}{\lvert}{\rvert}

\begin{document}
	\title{Spatial regularity for a class of degenerate Kolmogorov equations}
	
	\author{Francesca Anceschi}
	\address{Dipartimento di Matematica e Applicazioni
		"Renato Caccioppoli" -
		Università degli Studi di Napoli Federico II: 
		Via Cintia, Monte S. Angelo
		I-80126 Napoli, Italy}
	\email{francesca.anceschi@unina.it}
	
	\date{\today}

	\begin{abstract}
		\noindent
		We establish spatial a priori estimates for the solution $u$ to a class of
		dilation invariant Kolmogorov equation, where $u$ is assumed to only
		have a certain amount of regularity in the diffusion's directions, i.e. $x_{1},
		\ldots, x_{m_{0}}$. The result is that $u$ is also regular with respect to 
		the remaining directions. The approach we propose is based on the commutators identities 
		and allows us to obtain a Sobolev exponent that does not depend on the integrability assumption of the right-hand side.
		Lastly, we provide an alternative proof to that of Theorem 1.5 of \cite{Bouchut} for the optimal spatial regularity. 
				
		
	\end{abstract}
	
	\maketitle
	
\section{Introduction}
	The aim of this work is to prove spatial a priori estimates for solutions to a class of 
	strongly degenerate ultraparabolic equations
	\begin{align}\label{defL}
		\L u (x,t) &= Y u (x,t) - \sigma \Delta_{0} u \\ \nonumber
		&:= 
		\sum \limits_{i,j=1}^N b_{ij}x_j\partial_{x_i}u(x,t)-\partial_t u(x,t)  - \sigma \sum \limits_{i=1}^{m_{0}}
		\frac{\p^{2} u}{\p x_{i}^{2}}(x,t) = g(x,t),
	\end{align}
	where $(x,t)=(x_1,\ldots,x_N,t)\in \RR^{N+1}$, $1 \leq m_0 \leq N$ and $\sigma> 0$.
	Moreover, the matrix $B=(b_{ij})_{i,j=1,\ldots,N}$ satisfies the following structural assumption.

\medskip

\begin{itemize}
\item[\textbf{(H)}] For some basis on $\RR^N$, the matrix $B$ takes the following form
\begin{equation}
	\label{B}
	B =
	\begin{pmatrix}
		\OO   &   \OO      \\  
		B_1   &    \OO  
	\end{pmatrix}
\end{equation}
where $B_1$ is a $m_{1} \times m_{0}$ matrix of rank $m_1$, with 
\begin{equation} \label{defm}
	m_0 \ge m_1 \ge 1 \hspace{5mm} \text{and} \hspace{5mm} 
	m_0 + m_1 = N.
\end{equation} 
\end{itemize}
In particular, we provide spatial a priori estimates of the solution $u$ in terms of Sobolev spaces
when the right-hand side is in the $L^{2}$ space and $u$ is assumed to have a certain amount of
regularity in the first $m_{0}$ directions. The result is that $u$ is also regular in the remaining 
$N- m_{0}$ space variables and this is related to a commutator identity introduced by H\"ormander
for hypoelliptic operators. Results of this type where firstly proved in \cite{H} thanks to a representation formula via the 
fundamental solution. Later on, they were improved in \cite{RS} through a condition that considers brackets of the vector fields, but that leads to 
involved estimates. Recently, we recall the works \cite{Bouchut, BBB}, where the authors deal with the classical kinetic Fokker-Planck equation
and the class of homogeneous H\"ormander's operators that can be written as a sum of squares, i.e. $\sum X_i^2$, respectively.
Our aim is to extend the results of \cite{Bouchut} to the operator \eqref{defL}, via an elegant approach
based on the commutators identities that allows us to obtain a
Sobolev exponent that does not depend on the integrability assumption of the right-hand side.
Lastly, we provide an alternative proof to that of Theorem 1.5 of \cite{Bouchut} for the 
optimal spatial regularity.
The technique we consider is not based on any representation formula through the fundamental solution, 
or Fourier transform method, and allows us to obtain a first extension of 
the results in \cite{Bouchut} in the direction of the more general class of ultraparabolic operators 
\begin{align*}
	\K u (x,t) &:= Y u (x,t) + {\rm div} (A D u ) = g(x,t), \qquad (x,t) \in \RR^{N+1},
\end{align*}
known as Kolmogorov operators in divergence form, for which the matrix $B=(b_{ij})_{i,j=1,\ldots,N}$
that defines the transport operator $Y$ has the following form for some basis of $\RR^N$:
\begin{align}
	\label{Bk}
	B =
	\begin{pmatrix}
		\OO   &   \OO   & \ldots &    \OO   &   \OO   \\  
		B_1   &    \OO  & \ldots &    \OO    &   \OO  \\
		\OO    &    B_2  & \ldots &  \OO    &   \OO   \\
		\vdots & \vdots & \ddots & \vdots & \vdots \\
		\OO    &  \OO    &    \ldots & B_\k    & \OO
	\end{pmatrix}
\end{align}
where $B_i$ is a $m_i \times m_{i-1}$ matrix of rank $m_i$ for every $i=1, \ldots, \k$ and 
\begin{equation*}
	m_0 \ge m_1 \ge \ldots \ge m_\k \ge 1, \qquad \sum \limits_{i=0}^{\k} m_i = N.
\end{equation*}

\subsection{Motivation and background}
The Kolmogorov equation $\K u = g$ is a second order ultraparabolic PDE, that is strongly degenerate whenever $m_0 < N$ and 
arises in various different research fields, such as the financial market modeling \cite{PA, BPV, AMP}, or the
kinetic theory of gas \cite{IS-weak, AZ, DT}. It is known that the 
first order part of $\K$ may induce a strong regularizing property, and in fact was considered by H\"ormander in \cite{H} as a prototype for the 
family of hypoelliptic operators $\sum X_i^2 + X_0$. Indeed, under a particular structural assumption
for the matrix $B$ (see \eqref{Bk}), the operator $\K$ is homogeneous and hypoelliptic, see \cite{APsurvey, LP}. 
Namely, every distributional solution $u$ to $\L u = f$ defined in 
some open set $\Omega \subset \RR^{N+1}$ belongs to $C^\infty(\Omega)$ 
and it is a classical solution to $\K u = f$, 
whenever $f \in C^\infty(\Omega)$. Moreover, we remark that assumption \textbf{(H)} is implied by the condition introduced by H\"{o}rmander in \cite{H}:
\begin{equation*} 
{\rm rank\ Lie}\left(\p_1, \ldots, \p_{m_0},Y\right)(x,t) =N+1 \qquad \text{for every } \, (x,t) \in \RR^{N+1},
\end{equation*}
where ${\rm  Lie}\left(\p_1, \ldots, \p_{m_0},Y\right)$ denotes the Lie algebra generated by the first order differential operators $\left(\p_1, \ldots, \p_{m_0},Y\right)$ computed at $(x,t)$.

Nowadays, the research community is mainly interested in the study of the regularity theory for solutions to $\K u = g$, and in particular 
there are various open problems regarding the weak regularity theory, since the literature regarding the classical regularity theory is quite vast. 
For further information on this topic we refer to the survey paper \cite{APsurvey} and the references therein. 
On the other hand, the weak regularity theory has 
been widely investigated in the last decade and it is still evolving. 
Among the most recent results, we recall a series 
of paper devoted to the study of the Fokker-Planck operator \cite{GIMV, GI, GM, AM}, that is a particular case of the 
operator $\L$, i.e. 
\begin{equation*}
	\Delta_{v} u (v,x,t) + v \cdot \nabla_{v} u (v,x,t) - \p_{t} u (v,x,t) \qquad \text{for every } \, (v,x,t) \in \RR^{2d} \times \RR,
\end{equation*}
obtained by choosing $m_{0} = m_{1}=d$, $N = 2d$ and $B_{1} = \I_{d}$, where $\I_{d}$ is the identity matrix of order $d$.
Very recently, the authors of \cite{AR-harnack,WZ-preprint} extended the technique proposed in \cite{GI} for the Fokker-Planck equation and in \cite{IS-weak} for the Boltzmann equation to the case of weak solutions to $\K u = g$. 
This method is based on the combination of a weak Harnack inequality with a log-transform, and 
it is the only approach we can employ up to now when dealing with the case of weak solutions to $\K u = g$. 

Among the aforementioned open problems, we recall the extension of the De Giorgi
technique proposed in \cite{GIMV, GM} to weak solutions to the more general equation $\K u = g$. 
This method is based on the combination of an $L^{2}-L^{\infty}$ estimate ($1^{st}$ Lemma of De Giorgi) for weak solutions with a measure-to-pointwise 
result ($2^{nd}$ Lemma of De Giorgi), also known as intermediate value lemma for weak solutions. 
As far as we are concerned with the $1^{st}$ Lemma, this was proved in \cite{APR, AR-harnack} in a very general framework. 
On the other hand, in \cite{GIMV, GM} we find two different proofs for the intermediate value lemma for weak solutions to the classical
Fokker-Planck equation, but no generalization is available even for the case of \eqref{defL}.
In particular, in the paper \cite{GIMV} the proof of
 the $2^{nd}$ Lemma is based on global fractional estimates in Sobolev spaces proved by Bouchut in
\cite[Theorem 1.3]{Bouchut} via techniques from the velocity averaging theory, see \cite{DV} 
and the references therein. These techniques mainly consist of two different approaches - the Fourier transform and the commutators identities - 
and to our knowledge there is no available result of this type even for the class of \eqref{defL}.
Hence, the aim of this work is to provide a first extension of the results proved in \cite{Bouchut} 
to the equation \eqref{defL} following the commutators identities approach, since it allows us to treat more general operators without 
considering any representation formula via the fundamental solution, but with the only assumption of extra regularity in the 
degeneracy directions of $u$.

\subsection{Main results}
From now on, according to \eqref{defm} we split the coordinate $x\in\RR^N$ as
\begin{equation}\label{split.coord.RN}
	x=\big(x^{(0)},x^{(1)}\big), \qquad x^{(0)}\!\in\RR^{m_0}, \quad x^{(1)}\!\in\RR^{m_1}.
\end{equation}
Moreover, we denote by $\Delta_0$ the Laplacean with respect to the first $m_0$ directions
and the partial gradient in $\RR^{m_{i}}$ with respect to $x^{(i)}$ as
\begin{equation*}
	D_{i}= (\partial_{x^{(i)}_1},\ldots,\partial_{x^{(i)}_{m_i}}) 
	\qquad \text{for } \, i = 0, 1.
\end{equation*}
Lastly, we define for every $\beta \ge 1$ and for every $s \ge 0$
\begin{equation*}
	D^{\b}_{0} := (- \Delta_{m_{0}} u )^{\b/2} , \qquad 
	D^{s}_{1} := (- \Delta_{m_{1}} u )^{s/2}.
\end{equation*}
Now, we are in a position to state our main results, that respectively extend Theorem 1.6 and Theorem 1.5  of \cite{Bouchut}
to solutions to \eqref{defL}. Here, we do not explicitly express which notion of solution we consider, 
since our results apply to every solution of the equation \eqref{defL} that is at least in $L^{2} (\RR^{N} \times \RR)$. 
\begin{theorem} \label{thm-h1}
	Let us consider a solution $u \in L^2(\RR^{N} \times \RR)$ to
	\begin{equation*}
		Yu (x,t)  =g(x,t) \qquad  (x,t) \in \RR^{N} \times \RR,
	\end{equation*}
	 with $g \in L^{2}(\RR^{N} \times \RR )$ and 
	 \begin{equation*}
	 	D_0^{\max \{\b, \gamma\} } u \in L^{2}(\RR^{N} \times \RR ) \qquad \text{and}
		\qquad D_0^{\gamma } g \in L^{2}(\RR^{N} \times \RR )
	\end{equation*}
	with 
	\begin{equation*}
		\g \ge 0 \qquad \text{and} \qquad 0 \le 1 - \gamma \le \beta.
	\end{equation*}
	Then $D_1^{s} u \in L^{2}(\RR^{N+1})$, with 
	\begin{equation}
		\label{h2}
		s =  \frac{\b}{1 - \gamma + \beta } \qquad \text{and} \qquad s=1 \quad \text{if } \, 1 - \gamma = \beta = 0.
	\end{equation}
	Moreover, the following estimate holds
	\begin{align} \label{h3}
		 \| D_{1}^{s}u \|_{2} \le 
		C(N,B, \beta) \, \| D_{0}^{\beta} u \|^{1-s}_{2} \, \, \| D^{\gamma}_0 g \|_{2}^{s}.
	\end{align}
\end{theorem}
The proof we propose here is reminiscent of the hypoelliptic regularity proposed in \cite{H} and 
extends also \cite[Proposition 1.1]{Bouchut} ($\gamma = 0$), with the extra assumption that $\b \ge 1$.
This is due to the fact that the proof of \cite[Proposition 1.1]{Bouchut} is based on the application of the Fourier transform
to the equation $Y u = g$. 
From the weak regularity point of view, the most interesting case is given by the choice of $\gamma = 0$ and $\b = 1$, since 
$D_{0} u \in L^{2}( \RR^{N} \times \RR)$ and $g \in L^{2}( \RR^{N} \times \RR)$ are standard regularity assumptions when studying weak solutions to 
\eqref{defL}. In particular, this choice for the parameters leads to 
\begin{equation*}
	D_1^{\frac12} u \in L^{2}(\RR^{N+1}).
\end{equation*}
We point out that estimates of this type are very useful for the 
characterization of the space $\mathcal{W}$ of weak solutions introduced in \cite{AM, AR-harnack}, and
in order to complete the proof we only rely on the geometrical properties of the hypoelliptic structure of the transport operator. 

\medskip

When dealing with the ``complete'' evolution equation \eqref{defL}, it is quite clear how we can recover an estimate of 
$D_{0} u$ starting from the equation (see the Caccioppoli inequality proved in \cite[Theorem 3.1]{APR}). 
However, combining this information (that corresponds to the choice of $\b = 1$) with the previous theorem
does not provide us with the best Sobolev exponent for the spatial regularity, that in \cite{RS} is showed to be $2/3$ in the direction $x^{(1)}$. As we already 
pointed out at the beginning of this introduction, the method proposed in \cite{RS} is very involved. For this reason, in \cite[Theorem 1.5]{Bouchut}
the author proposed to show that a solution to \eqref{defL} is such that $\Delta_{m_{0}} u \in L^{2}$ and then apply Theorem \ref{thm-h1}, with $\b = 2$. 
Here, we state the analogous theorem for solutions to \eqref{defL}, for whose proof we refer to \cite{Bouchut}. 
\begin{theorem} \label{thm-h3}
	Let us consider a solution $u \in L^2(\RR^{N} \times \RR)$ to 
	\begin{equation*}
		Yu (x,t) - \sigma \Delta_{m_0} u (x,t) =g(x,t) \qquad  (x,t) \in \RR^{N} \times \RR,
	\end{equation*}
	where $\sigma > 0 $ a positive constant. If $g \in L^2(\RR^N \times \RR)$ and $D_{0}^{\beta} u \in L^2(\RR^N \times \RR)$,
	with $\beta \ge 1$. Then, $Y u$ and $\Delta_{m_0} u$ are in $L^2(\RR^N \times \RR)$ and
	\begin{equation*}
		\| Y u \|_2 + \s \| \Delta_{m_0} u \|_2 \le C(N, B) \| g \|_2. 
	\end{equation*}	
	Moreover, $D_1^{\frac23} u \in L^2(\RR^N \times \RR)$ and
	\begin{equation*}
		\| D_1^{\frac23} u \|_2 \le \frac{C(N,B)}{\s^{1/3}} \| g\|_2.
	\end{equation*}
\end{theorem} 

Lastly, we propose an alternative proof of the maximal regularity carried out by contradiction, 
whose main advantage is it directly extends to the case of $\K u =g$ once the corresponding spatial estimates are proved. 
\begin{theorem} \label{thm-h2}
	Let us consider a solution $u \in L^2(\RR^{N} \times \RR)$ to 
	\begin{equation*}
		Yu (x,t) - \sigma \Delta_{m_0} u (x,t) =g(x,t) \qquad  (x,t) \in \RR^{N} \times \RR,
	\end{equation*}
	with $g \in L^2(\RR^N \times \RR)$ and $\sigma > 0 $ a positive constant such that $\sigma \ne 1$. 
	Then, both $Y u$ and $\Delta_{m_0} u$ are in $L^2(\RR^N \times \RR)$ and the estimate \eqref{h3}
	holds true with the choice of the parameter $\b=2$.	
\end{theorem} 
On one hand the contradiction argument forces us to require $\s \ne 1$, on the other hand it allows us to get rid of the 
assumption of $D_{0}^{\b} u \in L^{2}(\RR^{N} \times \RR)$ that we have in Theorem \ref{thm-h3}.
Also, we believe it is possible to extend this approach to the more general case of the Kolmogorov equation in divergence (and non-divergence form) with measurable coefficients, for which it is not possible to consider an approach in the spirit of \cite{H} since the existence of a fundamental solution associated to it has not been proved yet. 

\medskip

This paper is organized as follows. 
Section \ref{proof1} and Section \ref{proof2} are respectively devoted to the 
proof of Theorem \ref{thm-h1} and Theorem \ref{thm-h2}. In Section \ref{ex} 
we state our conclusions and we introduce the model operator of the Kolmogorov equation with two groups of spatial 
variables.

\setcounter{equation}{0}\setcounter{theorem}{0}
\section{Proof of Theorem \ref{thm-h1}} \label{proof1}
By suitable smoothing and cut-off, we can reduce ourselves to the case where $u,g$ are smooth 
functions, i.e. $u,g \in C^{\infty}_{c}(\RR^{N+1})$. Following the notation of \cite{H}, 
from now on we denote the $L^{2}$ pairing by
\begin{equation} \label{bracket}
	\langle f, g \rangle = \int \limits_{\RR^{N+1}} fg \, dx \, dt . 
\end{equation}
By assumption \textbf{(H)}, the block $B_1$ of the matrix $B$ is of maximum rank 
$m_{1}$, since by \eqref{defm} we have that $m_1 \le m_{0}$. Without loss of
generality, we may assume the block $B_1$ is in the canonical form with $m_1$ pivots, i.e.
\begin{align}
	\label{pivot}
	B_1 = \begin{pmatrix}
			b^1_{11} & b^1_{12} & \ldots & b^1_{1 m_{1}} & \ldots & b^1_{1 m_{0}}\\
			0 & b^1_{22} & \ldots & b^1_{2 m_{1}} & \ldots & b^i_{2 m_{0}} \\
			\vdots & \vdots & \ddots & \vdots & \ddots & \vdots \\
			0  & \ldots  & \ldots & b^1_{m_{1} m_{1}} & \ldots & b^1_{m_{1} m_{0}} 
	\end{pmatrix}
\end{align}
The proof of Theorem \ref{thm-h1} is based on the representation 
of the derivative along a certain direction $x^{(1)}_{j}$, through a commutators identity between 
$\p_{x^{(0)}_{j}}$ and $Y$.
 In particular, thanks to the explicit expression of $Y$, we have that
\begin{align} \label{comm-id}
	[\p_{x^{(0)}_{1}} , Y ] u = \p_{x^{(0)}_{1}}  Y u - Y \p_{x^{(0)}_{1}}  u &
	= b_{1 1}^1 \p_{x^{(1)}_{1}} u \\ \nonumber
	[\p_{x^{(0)}_{2}} , Y ] u = \p_{x^{(0)}_{2}}  Y u - Y \p_{x^{(0)}_{2}}  u &=  
	b_{1 2}^1 \p_{x^{(1)}_{1}} u  + b_{2 2}^1 \p_{x^{(1)}_{2}} u \\ \nonumber
	&\, \, \vdots \\ \nonumber
	[\p_{x^{(0)}_{m_{1}}} , Y ] u = \p_{x^{(0)}_{m_{1}}}  Y u - Y \p_{x^{(0)}_{m_{1}}}  u &=  
	\sum \limits_{j=1}^{m_{1}} b^{1}_{j m_1} \p_{x^{(1)}_j} u.
\end{align}
First of all, this procedure is well posed since $m_{0} \ge m_{1}$ and the matrix $B_{1}$ is of maximum rank
$m_1$. Hence, thanks to the first commutators identity of \eqref{comm-id}, we have 
\begin{align*}
	\p_{x^{(1)}_{1}} u  = \frac{1}{b_{1 1}^{(1)}} [\p_{x^{(0)}_{1}} , Y ] u 
\end{align*}
Then, there exists two constants such that 
\begin{equation} \label{es1}
		\p_{x^{(1)}_{2}} u = c_1 \p_{x^{(1)}_{1}} u + c_2 [\p_{x^{(0)}_{2}} , Y ] u .
\end{equation}
Generally speaking, for a fixed $i =1, \ldots, m_1$, 
there exist $c_1, \ldots, c_{i}$ constants such that
\begin{equation} \label{es2}
\p_{x^{(1)}_i} u =  \sum \limits_{j=1}^{i-1} c_j \p_{x^{(1)}_{j}} u + c_i [\p_{x^{(0)}_{i}} , Y ] u.
\end{equation}

Now, the idea behind the proof of this theorem is to consider the commutators 
identities above and to estimate the $L^{2}$ bracket of each of them against a suitable derivative. 
Let us begin with the estimate of the derivative $\p_{x^{(1)}_1} u$:
\begin{align} \label{es3}
	\| D_1^{-(1-s)} \p_{x^{(1)}_1} u \|_{2}^{2} &= 
	\langle D_1^{-2(1-s)} \p_{x^{(1)}_1} \overline u, \p_{x^{(1)}_1} u \rangle 
	\\ \nonumber
	 &= \langle D_1^{-2(1-s)} \p_{x^{(1)}_1} \overline u, \p_{x^{(0)}_1} Y u - Y  \p_{x^{(0)}_1} u \rangle 
	 \\ \nonumber
	&= \langle D_1^{-2(1-s)} \p_{x^{(1)}_1} \overline u, \p_{x^{(0)}_1}g  - Y  \p_{x^{(0)}_1} u \rangle 
	 \\ \nonumber
	&= -  \langle  \p_{x^0_1}  D_1^{-2(1-s)}  \p_{x^{(1)}_1}\overline u,  g \rangle
	 + \langle Y D_1^{-2(1-s)}  \p_{x^{(1)}_1}  \overline u, \p_{x^{(0)}_1} u \rangle \\ \nonumber
	&= - 2 Re \langle D_1^{-2(1-s)}  \p_{x^{(1)}_1} \p_{x^{(0)}_1} \overline u,  g \rangle \\ \nonumber
	&\le C(N,B)  \, \| D_1^{1-2(1-s)}  D_0 u \|_{2} \, \| g \|_{2}.
\end{align}
Moreover, thanks to \eqref{es1} we also have 
\begin{align*}
	\| D_1^{-(1-s)} \p_{x^{(1)}_2} u \|_{2}^{2} &= 
	\langle D_1^{-2(1-s)} \p_{x^{(1)}_2} \overline u, \p_{x^{(1)}_2} u \rangle 
	\\ \nonumber
	&= c_1   \langle D_1^{-2(1-s)} \p_{x^{(1)}_2} \overline u,  \p_{x^{(0)}_1} Y u - Y  \p_{x^{(0)}_1} u \rangle \\ \nonumber
	&\qquad \qquad +
	c_2 \langle D_1^{-2(1-s)} \p_{x^{(1)}_2} \overline u, \p_{x^{(0)}_2} Y u - Y  \p_{x^{(0)}_2} u \rangle \\ \nonumber
	&\le C(N, B)   \, \| D_1^{1-2(1-s)}  D_0 u \|_{2} \, \| g \|_{2}.
\end{align*}
Taking into consideration \eqref{es2}, the above argument holds true also for every other direction $x^{(1)}_i$, with $i=2, \ldots, m_1$. For this reason, 
from now on we restrict ourselves to the case where $B_1$ is of the form \eqref{pivot}, where every component is zero except $b^1_{ii}=1$, with 
$i=1, \ldots, m_1$, i.e. 
\begin{align*}
	B_1 = \begin{pmatrix}
			1 & 0 & \ldots & 0 & \ldots & 0\\
			0 & 1& \ldots & 0 & \ldots & 0 \\
			\vdots & \vdots & \ddots & \vdots & \ddots & \vdots \\
			0  & \ldots  & \ldots & 1 & \ldots & 0
	\end{pmatrix},
\end{align*}
and we keep in mind that, from now, every constant we consider depends on the matrix $B$.
By noticing that $\frac12 \le s < 1$ and that \eqref{es3} can be computed for every $\p_{x^{(1)}_i} u$, with 
$i=1, \ldots, m_1$, we have that:
\begin{align*}
	\| D_1^{s} u \|_{2}^{2}  \le 
	C(N,B)  \, \| D_1^{1-2(1-s)}  D_0^{1-\gamma} u \|_{2} \, \| D_0^{\gamma} g \|_{2}.
\end{align*}
Then we define $\theta = 2 - 1/s \in [0,1]$ and notice that $1-2(1-s) = \theta s$ and $1-\gamma = (1-\theta) \beta$.
Then if we apply the H\"older inequality in the Fourier variables we consider to define the fractional derivative in $x^{(1)}$ (see equation (3.6) of \cite{Bouchut}), we get 
\begin{align}
	\label{ref-c}
	\| D_1^{s} u \|_{2}^{2} &\le
	C(N,B) \, \| D_1^{s} u \|_{2}^{\theta} \, \, \| D_0^{\b} u \|^{1-\theta}_{2} \, \, \| D_0^{\gamma} g \|_{2}.
\end{align}
Lastly, by collecting together all of the corresponding terms, we obtain
\begin{align*}
	\|  D_1^{s}u \|_{2}^{2-\theta} &\le
	C(N,B) \, \|  D_0^{\b} u \|^{1-\theta}_{2} \, \, \| D_0^{\gamma} g \|_{2}.
\end{align*}	
Hence, the proof is complete.
$\hfill \square$

\setcounter{equation}{0}\setcounter{theorem}{0}
\section{Proof of Theorem \ref{thm-h2}}
\label{proof2}
We proceed by contradiction. Thus, we assume $Yu \notin L^2(\RR^{N+1})$ and $\Delta_{m_0} u \in L^2(\RR^{N+1})$.
By computing the $L^2$ norm on both sides of \eqref{defL} and by applying the inverse triangle inequality we have that 
\begin{equation*}
	\abs{ \| Y u \|_2 - \sigma \| \Delta_{m_0} u\|_2 } \le \| Yu - \s  \Delta_{m_0} u \|_2 = \| g \|_2.
\end{equation*}
Thus, by considering that $\| Yu \|_2$ is unbounded, the above inequality leads us to a contradiction:
\begin{equation*}
	+ \infty \le \| g \|_2 .
\end{equation*}
The same reasoning applies if we assume $\Delta_{m_0} u \notin L^2(\RR^{N+1})$ and $Y u \in L^2(\RR^{N+1})$.
Lastly, if we assume both $\Delta_{m_0} u, Yu \notin L^2(\RR^{N+1})$, the only interesting case is when their asymptotic
behavior is comparable. Still, since by assumption $\s > 0$ and $\s \ne 1$, what we obtain is a contradiction again. Indeed, if we denote by
$\o = \| \Delta_{m_0} u \|_2 = \|Yu \|_2 $ and we consider the reverse triangle inequality once more, we obtain
\begin{equation*}
	+ \infty = \lim \limits_{\o \to + \infty} \abs{1 - \sigma } \o  \le  \| g \|_2.
\end{equation*}
$\hfill \square$

\setcounter{equation}{0}\setcounter{theorem}{0}
\section{Conclusions}
\label{ex}
We have fully addressed the first issue we encounter when extending the results of \cite{Bouchut} 
to the Kolmogorov equation $\K u = g$, that is represented by the choice of the matrix $B$ as in \textbf{(H)}, i.e.
with $m_0 \ge m_1$, $m_0$ possibly different from $m_1$, and $B_1$ a generic matrix of order $m_1 \times m_0$, 
not necessarily equal to the identity. 

The next issue we have to face is to consider a prototype equation for the family of $\K u = g$ that introduces 
the mixing structure of the transport operator $Y$. In particular, we are interested in the preliminary study of the following toy model equation, 
also known as Kolmogorov equation with two groups of spatial variables:
	\begin{align}\label{toy}
		\widetilde Y u (x,t) - \sigma \Delta_{0} u(x,t) 
		:= & \sum \limits_{j=1}^d x^{(1)}_{j} \partial_{x^{(0)}_{j}} u(x,t) +
		\sum \limits_{j=1}^d x^{(2)}_{j} \partial_{x^{(1)}_{j}} u(x,t) \\ \nonumber
		&\qquad - \partial_t u(x,t)  
		- \sigma \sum \limits_{i=1}^{m_{0}}
		\frac{\p^{2} u}{\p x_{i}^{2}}(x,t) = g(x,t),
	\end{align}
	where $(x,t)=(x^{(0)},x^{(1)}, x^{(2)}, t)\in \RR^{3d+1}$ and $x^{(i)} \in \RR^{d}$ for every $i=1,2,3$.
	The toy model \eqref{toy} is obtained by considering the following choice for the 
	matrix $B$
	\begin{equation*}
	\label{Bex}
	B =
	\begin{pmatrix}
		\OO_{d}   &   \OO_{d}  &  \OO_{d}    \\  
		B_1   &    \OO_{d}  &  \OO_{d}   \\
		\OO_{d}    &    B_2    &  \OO_{d} 
	\end{pmatrix},
\end{equation*}	
where for $i=1,2$ every $B_{i}= \I_d$.
This toy model can be seen as a first extension of the Fokker-Planck equation where we have one velocity direction $x^{(0)}$ and 
two spatial directions $x^{(1)}$ and $x^{(2)}$ and, from a geometrical point of view, presents all of the difficulties we encounter when studying the more 
general Kolmogorov equation $\K u = g$. Indeed, even if the matrix $B$ is simply made of $B_{1}$, $B_{2}$ and two non-zero blocks, it is immediately clear 
how the structure of the transport operator 
\begin{equation*}
		\widetilde Y u = \sum \limits_{j=1}^d x^{(1)}_{j} \partial_{x^{(0)}_{j}} u(x^{(1)}, x^{(2)},t) +
		\sum \limits_{j=1}^d x^{(2)}_{j} \partial_{x^{(1)}_{j}} u(x^{(1)}, x^{(2)},t)-\partial_t u(x^{(1)}, x^{(2)},t) 
\end{equation*}
mixes the spatial direction $x^{(1)}$ with the velocity direction $x^{(0)}$, but also the spatial direction 
$x^{(2)}$ with the velocity direction $x^{(1)}$. On one hand, this means it is not straightforward to adapt the kinetic theory's techniques 
to \eqref{toy}. On the other hand, it is clear that it presents all the difficulties one has to deal with in the more general case of $\K u = g$.
For a recent bibliography on the study of this particular model we refer to \cite{IM-toy1, IM-toy2} and the reference therein.

In the future, our aim is to prove a spatial regularity result in the spirit of Theorem \ref{thm-h1}
for the toy model \eqref{toy}, and then to consider it to prove analogous results to that of Theorem \ref{thm-h3} and 
Theorem \ref{thm-h2}. A first attempt to the proof of this result would be to consider the technique presented in Section \ref{proof1},
with $s = \frac{\b}{\b+1}$ that plays the same role as $\b$ in the estimate \eqref{ref-c} when computing the iteration for the second space variable 
$x^{(2)}$. Still, this is not well posed since the derivative we would obtain in the last line of \eqref{es3} would have a negative exponent. 
On the other hand, we conjecture it is possible to extend the Fourier technique considered in \cite[Theorem 1.3]{Bouchut} to the toy model \eqref{toy}.

\section*{Aknowledgments} 
The author is funded by the research grant PRIN2017 2017AYM8XW ``Nonlinear Differential problems via variational, topological and set valued methods''.

\end{document}